\documentclass[pdflatex,sn-mathphys-num]{sn-jnl}

\usepackage{graphicx}%
\usepackage{multirow}%
\usepackage{amsmath,amssymb,amsfonts}%
\usepackage{amsthm}%
\usepackage{mathrsfs}%
\usepackage[title]{appendix}%
\usepackage{xcolor}%
\usepackage{textcomp}%
\usepackage{manyfoot}%
\usepackage{booktabs}%
\usepackage{algorithmic}%
\usepackage{algorithm}%
\usepackage{listings}%

\theoremstyle{thmstyleone}%
\newtheorem{theorem}{Theorem}
\newtheorem{definition}{Definition}
\theoremstyle{thmstyletwo}%

\theoremstyle{thmstylethree}%

\newcommand{\x}{\mathbf{x}}
\newcommand{\pos}{\mathbf{r}}
\newcommand{\vel}{\mathbf{v}}
\newcommand{\acc}{\mathbf{a}}
\newcommand{\pv}{\mathbf{p}}
\newcommand{\uc}{\mathbf{u}}
\newcommand{\R}{\mathbb{R}}

\newcommand{\bPhi}{\boldsymbol{\Phi}}
\newcommand{\bGamma}{\boldsymbol{\Gamma}}
\newcommand{\bPsi}{\boldsymbol{\Psi}}

\begin{document}

\title[Multi Impulse Low Earth Orbit Maneuver Synthesis Through Thrust Measure Primer Vector Conditions and Piecewise Multiple Shooting]{Multi Impulse Low Earth Orbit Maneuver Synthesis Through Thrust Measure Primer Vector Conditions and Piecewise Multiple Shooting}


\author*[1]{\fnm{Pedro} \sur{Kuntz Puglia}}\email{pepuglia@gmail.com}

\author[1]{\fnm{Willer} \sur{Gomes dos Santos}}\email{willer@ita.br}
\equalcont{These authors contributed equally to this work.}

\author[2]{\fnm{Emilien} \sur{Flayac}}\email{Emilien.FLAYAC@isae-supaero.fr}
\equalcont{These authors contributed equally to this work.}


\affil[1]{\orgdiv{IEA},  \orgname{ITA}, \orgaddress{\street{Praça Marechal Eduardo Gomes}, \city{São José dos Campos}, \postcode{12228-900}, \state{São Paulo}, \country{Brazil}}}

\affil[2]{\orgdiv{Fédération ENAC ISAE-SUPAERO ONERA}, \orgname{Université de Toulouse}, 
\orgaddress{\street{10 Av. Marc Pélegrin}, \city{Toulouse}, \postcode{31055}, \country{France}}}


\abstract{This work sets out to apply primer vector theory to optimal impulse maneuvers in Low Earth Orbit (LEO) and adapt it to the types of perturbations encountered in this environment, which is not readily available in the literature. A review of the theory of optimal control and orbital maneuvering is made and in particular, primer vector theory is laid out in detail based on the extension of the thrust control to a measure, including its generalization to conservative and non-conservative perturbation models. An impulsive multiple shooting optimization scheme in Cartesian coordinates is presented, through a piecewise approach to the problem. Then, some maneuver scenarios with known solutions under the Keplerian model are solved under perturbed orbital dynamics with the help of the primer vector, and the resulting trajectories are compared. The perturbation models include a J2 model, representing the class of conservative perturbations by modeling Earth's nonspherical gravity field, and a J2+Drag model, representing the class of non-conservative perturbations by including the effects of atmospheric drag in LEO\@. For each model, some valid methods of primer vector calculation have been tried and validated between each other, and the primer vector is proven to be a useful tool in reducing the cost of orbital maneuvers.}

\keywords{Orbital Maneuvering, Optimal Control, Impulsive thrust, LEO}



\maketitle

\section{Introduction}\label{sec1}

The issue of orbital maneuvering goes back to the early days of space exploration. The prevalent objective of minimizing fuel consumption is deceptively simple, and leads to challenging optimization problems whose solution is heavily dependent on problem details. Many common variants exist, differing in the propulsive model, the maneuver time, fixed or free, and initial and final conditions, as well as the orbital model, such as two- or multi-body problems.

Low Earth Orbit (LEO) maneuvers are of particular interest to Brazil's current space program. The ITASAT-2 project proposes the formation flying of 3 16U Cubesats in LEO that must keep precise positions through the usage of chemical thrusters~\cite{shibuya_sato_itasat-2_2023}. While previous work considers linear closed-loop relative position control in a small deviations context~\cite{franco_itasat-2_2020}, this work deals with the complimentary problem of nonlinear open loop absolute position trajectory generation. The LEO environment imposes dynamical perturbations that change the behavior of trajectories. In particular, Earth's non-spherical gravity field makes trajectories nonplanar and atmospheric drag makes orbits decay~\cite{curtis_orbital_2020}. Particularly important for this work is the fact that the LEO environment. Existing orbital maneuvering works rarely include these perturbations in the dynamical model.

A common approach for modeling propulsive systems is that of impulsive thrust. Physically, it is an approximation for satellites with high acceleration capabilities, such as those with chemical rocket engines~\cite{conway_spacecraft_2010}. Mathematically, removing the bound on the maximum allowable thrust results in  unbounded delta-V cost when considering finite thrust trajectories; treating the problem in a rigorous way requires the extension of the thrust control from a bounded function to a measure~\cite{arutyunov_optimal_2019}. This thrust measure approach relies on modern mathematical results and is not studied in the orbital maneuvering literature, and is one of the main contributions of this work.

Different boundary condition choices lead to very different solution types. Some results consider initial and final states belonging to manifolds. Notably, the Hohmann and bielliptic transfers, known for their analytical solutions~\cite{curtis_orbital_2020}, are free-time manifold to manifold solutions. Other works impose a fixed maneuver time and fixed boundary states, in an approach that is closer to a rendez-vous problem and is often used in the literature~\cite{prussing_optimal_1986, luo_interactive_2010}. 

The central question when dealing with impulsive maneuvers is how many impulses are necessary~\cite{conway_spacecraft_2010}. This is a discrete parameter which cannot be optimized numerically in the same way as the others. Some techniques for approaching this issue rely on homotopic continuation of bounded thrust maneuvers, where the maximum thrust is progressively increased~\cite{taheri_how_2020, arya_generation_2023}. These approaches circumvent the issue of finding the optimal number of impulses through this finite thrust parameterization, at the cost of obtaining an approximate solution.

Another common approach is using the \textit{primer vector}~\cite{lawden_optimal_1963}, a name given to the costate that arises in the indirect analysis of orbital maneuvering through the Hamiltonian formalism and Pontryagin's Maximum principle~\cite{conway_spacecraft_2010}. In the pioneering work that introduced the primer vector, Lawden coined the term by analogy with ignition primers, since the primer vector trajectory can indicate whether a maneuver is locally optimal or not~\cite{lawden_optimal_1963}. Since then, it has been extended in multiple ways. 

The primer vector optimality conditions have been extended to include diagnostic information about suboptimal trajectories, that is, how they could be modified to reduce cost, through calculus of variations analyses~\cite{lion_primer_1968}. These conditions were later used as the basis for an iterative algorithm that adds impulses and coasting arcs to a maneuver until the optimality conditions are satisfied~\cite{luo_interactive_2010, jezewski_efficient_1968}.

The primer vector method has been applied to a variety of problems. Many early and modern works apply it to the Keplerian (two body) model~\cite{lawden_optimal_1963, jezewski_efficient_1968, prussing_optimal_1986, luo_interactive_2010, jamison_analytical_2010, rebelo_optimizing_2024, sarli_orbit_2015, lion_primer_1968}, and it is widely explored in derived models such as multi-body problems~\cite{bucchioni_optimal_2023, bokelmann_optimization_2020, glandorf_primer_1970, bell1995primer} and linearized rendez-vous systems~\cite{serra_analytical_2018, zheng_optimal_2024, aubin_optimization_2011, kara-zaitri_polynomial_2010}. Lacking in the literature are applications to Low Earth Orbit, considering the characteristic orbit perturbations that affect this type of orbit, namely J2 effects and atmospheric drag.

Through analyzing the non-conservative dynamical model of a LEO orbit including atmospheric drag, a distinction can be introduced between the state transition matrix (STM) and the primer vector transition matrix (PVTM), often taken to be equal without further question~\cite{conway_spacecraft_2010, bokelmann_optimization_2020} since to this day the primer vector literature has, to the authors' knowledge, only dealt with conservative orbital models. The objective of this work is therefore to extend primer vector theory to the LEO environment, and illustrate how the orbital perturbations encountered there can change orbital maneuvers.

This paper is organized as follows. Section~\ref{sec:orb_models} introduced the orbital models considered, including the orbital maneuvering model and its general optimal control problem. Section~\ref{sec:piecewise} introduces a piecewise approach for the impulsive system, conducive to numerical implementation. Section~\ref{sec:pv_intro} introduces the impulsive optimality conditions which lead to a more rigorous derivation of primer vector theory, as well as the extension thereof to non-conservative models, and is the main contribution of this work. Section~\ref{sec:impl} details the optimization problem implemented numerically. Finally, in Section~\ref{sec:res} some numerical results are computed to illustrate the exposed theory, and Section~\ref{sec:conc} concludes the paper. 

\section{Orbital models}\label{sec:orb_models}

This Section shall explore different orbital models and introduce the orbital maneuvering optimal control problem. All orbital models include the central body's gravitational force, and differ on how orbital perturbations are treated. Keplerian motion is the motion of a satellite moving under spherically symmetric gravitational attraction~\cite{chobotov_orbital_2002}. In LEO, the most important perturbations are the Earth's oblateness effects, in particular the first zonal harmonic coefficient $J_2$, and atmospheric drag~\cite{curtis_orbital_2020}. Let \(\pos \in \R^3\) and \(\vel \in \R^3\) be the position and velocity of the satellite in Cartesian coordinates, in the Earth Centric Inertial Frame, and $\x = \begin{bmatrix}
    \pos^T & \vel^T
\end{bmatrix}^T$ be the full state vector. When dealing with any orbital model, the dynamics will be denoted
\begin{equation}
    \dot{\x} = \mathbf{f}(\x) = \begin{bmatrix}
        \dot{\pos} \\ \dot{\vel}
    \end{bmatrix} = \begin{bmatrix}
        \vel \\ \acc(\x)
    \end{bmatrix},
\end{equation}

where $\mathbf{f}: \R^6 \rightarrow \R^6$ is the dynamical function and $\acc: \R^6 \rightarrow \R^3$ is the satellite's acceleration. Now, the classes of conservative and non-conservative models are defined and exemplified.


\begin{definition}[conservative and non conservative orbital models]\label{def:cons_orb_model}
    An orbital model $\acc$ is said to be \emph{conservative} if there exists  a smooth gravitational potential field $V \in \mathcal{C}^2(\R^3 / \{\mathbf{0}_3\})$ such that:
\begin{align*}
    \acc(\x) = \acc(\pos) = \nabla V(\pos).
\end{align*}
Otherwise, it is said to be non conservative and is denoted by $\acc(\mathbf{x})=\acc(\mathbf{r},\mathbf{v})$.
\end{definition}

The Keplerian model described by the dynamics
\begin{equation}\label{eq:kepler_dyn}
    \ddot{\mathbf{r}} = -\frac{\mu}{\lVert \mathbf{r} \rVert^3} \mathbf{r},
\end{equation}
where \(\mu > 0\) is the gravitational parameter of the central body, is the simplest conservative model. The $J_2$ perturbation adds a conservative acceleration term \(\acc_{J2}(\pos)\) to the Keplerian model, given by~\cite{curtis_orbital_2020}
\begin{equation}\label{eq:acc_j2}
    \acc_{J2}(\pos) = \frac{3 J_2 \mu R^2}{2 r^4} \begin{bmatrix}
        \frac{x}{r} \left(5 \frac{z^2}{r^2} - 1\right) \\
        \frac{y}{r} \left(5 \frac{z^2}{r^2} - 1\right) \\
        \frac{z}{r} \left(5 \frac{z^2}{r^2} - 3\right)
    \end{bmatrix},
\end{equation}
where \(J_2 \in \R\) is the Earth's first zonal oblateness coefficient and \(R > 0\), is the equatorial radius~\cite{curtis_orbital_2020}. 


The drag perturbation is the most common non-conservative force in orbital mechanics, and is important in LEO. It relies on an atmospheric density model, \(\rho(r)\). A simple model is given by the US Standard Atmosphere 1976~\cite{curtis_orbital_2020}, which gives a piecewise exponential interpolation of empirical data. The atmosphere rotates with the Earth with angular velocity \(\boldsymbol{\omega}_E > 0\). The satellite is parameterized by a reference surface \(S > 0\), drag coefficient \(C_D > 0\) and mass \(m > 0\). The relative velocity of the satellite with respect to the atmosphere is given by
\begin{equation}
    \vel_{r} = \vel - (\boldsymbol{\omega}_E) \times \pos,
\end{equation}
and the acceleration due to drag \(\acc_D\) is finally given by~\cite{curtis_orbital_2020}
\begin{equation}\label{eq:acc_drag}
    \acc_D(\pos, \vel) = - \frac{1}{2} \rho(r) \frac{S C_D}{m} \lVert \vel_r \rVert \vel_r,
\end{equation}
which can be added as a perturbing acceleration to a gravitational model. Since drag is typically smaller than $J_2$ effects in LEO, drag should only be modeled if $J_2$ is included.

The distinction between conservative and non-conservative orbital models will be important for the primer vector discussion in Section~\ref{sec:pv_alg}. Table~\ref{tab:orb_models} summarizes these modeling considerations.

\begin{table}[htbp]
    \centering
    \begin{tabular}{cccc}\toprule
        Model & Model Class & Dynamics \\ \midrule
        Keplerian & Conservative & \(\ddot{\pos} = -\frac{\mu}{r^3} \pos\) \\
        J2 & Conservative & \(\ddot{\pos} = -\frac{\mu}{r^3} \pos + \acc_{J2}(\pos)\) \\
        J2+Drag & Non-conservative & \(\ddot{\pos} = -\frac{\mu}{r^3} \pos + \acc_{J2}(\pos) + \acc_{D}(\pos, \vel)\) \\ \bottomrule
    \end{tabular}
    \caption{Summary of orbital models, classes and dynamics considered.}
    \label{tab:orb_models}
\end{table}

For the \(n\)-dimensional dynamical system \(\dot{\x} = \mathbf{f}(\x)\), define the state transition matrix (STM) \(\bPhi_\delta(t) \in \R^{n \times n}\) in terms of the solution \(\x(t)\) and the initial condition \(\x_0\)
\begin{equation}\label{eq:stm_def}
    \mathbf{\Phi}_\delta(t) = \frac{\partial \x(t)}{\partial \x_0}.
\end{equation}

This matrix is subject to the initial value problem~\cite{montenbruck_satellite_2002}
\begin{equation}\label{eq:stm_ivp}
    \begin{cases}
        \dot{\mathbf{\Phi}}_\delta(t) = \frac{\partial \mathbf{f}(\x(t))}{\partial \x(t)} \mathbf{\Phi}_\delta(t) \\
        \mathbf{\Phi}_\delta(0) = \mathbf{I}_n
    \end{cases}
\end{equation}
with \(n^2\) variables and initial conditions. 

For the Keplerian model, an analytical expression for the STM exists, given for instance in \cite{montenbruck_satellite_2002} and \cite{glandorf_lagrange_1969}. It introduces an intermediary matrix \(\mathbf{P}(t)\), with analytical inverse \(\mathbf{P}^{-1}(t)\) such that
\begin{equation}\label{eq:stm_glandorf}
    \bPhi_\delta(t) = \mathbf{P}(t) \mathbf{P}^{-1}(0).
\end{equation}

For perturbed models, it is necessary to numerically integrate the initial value problem in Equation~\eqref{eq:stm_ivp} alongside the states or, if a convenient way of computing derivatives of ordinary differential equation (ODE) solutions is available, Equation~\eqref{eq:stm_def} can be used directly.

Note that any linear ODE with state \(\mathbf{z}(t) \in \R^n\) in the form
\begin{equation}
    \dot{\mathbf{z}}(t) = \frac{\partial \mathbf{f}(\x(t))}{\partial \x(t)} \mathbf{z}(t)
\end{equation}
has a solution in terms of its initial condition \(\mathbf{z}_0\) given by~\cite{montenbruck_satellite_2002}
\begin{equation}
        \mathbf{z}(t) = \bPhi_\delta(t) \mathbf{z}_0.
\end{equation}

\subsection{Orbital Maneuvering}\label{sec:orb_man}

In this section, propulsive control variables will be introduced in the dynamical model to allow for orbital maneuvers. The generic orbital model $\dot{\x} = \mathbf{f}(\x)$ will be used throughout the section, as the considerations in this section apply to all models. 

Several variations of the orbital maneuvering problem exist. Free time problems with boundary states belonging to fixed orbits lead to classical analytical solutions such as the Hohmann and bielliptic transfers~\cite{curtis_orbital_2020}. In numerical works, a fixed maneuver time is commonly chosen~\cite{jezewski_efficient_1968, serra_analytical_2018, arya_generation_2023, luo_interactive_2010, zheng_optimal_2024, prussing_optimal_1986, bucchioni_optimal_2023, rebelo_optimizing_2024, lion_primer_1968, kara-zaitri_polynomial_2010} to help convergence~\cite{bokelmann_optimization_2020}, and fixed boundary states may also be fixed. In this work, a fixed maneuver time $t_f > 0$ is assumed, as well as fixed Cartesian boundary states $\x_o \in \R^6$ and $\x_f \in \R^6$ for the initial and final conditions respectively.

The full continuous time impulsive maneuvering optimization problem is now stated in a standard form of the optimal impulsive control literature. Let
\begin{itemize}
    \item \(\Gamma\) be a scalar non-negative Radon measure on \([0, t_f]\) representing the acceleration magnitude;
    \item \(\hat{\uc}: [0, t_f] \rightarrow S^2\) be the trajectory of thrust directions;
    \item  \(\x: [0, t_f] \rightarrow \R^6\) be the bounded variation orbital state trajectory;
    \item and \(\Delta v: [0, t_f] \rightarrow \R\) be an extra bounded variation state representing the accumulated delta-V cost of the maneuver.
\end{itemize}

The acceleration vector measure $\bGamma$ is given by $d\bGamma = \hat{\uc} d\Gamma$. Then, the goal is to find the optimal cost \(\Delta v^*\):

\begin{align}\label{eq:measure_opt_prob}
    P_\Gamma = \begin{cases}
    \begin{tabular}{ccl}
     $\Delta v^\star = \min$ & $\Delta v(t_f)$\\
         $\Gamma, \hat{\uc}, \x, \Delta v$ &  \\
          subject to: &$\x(0) = \x_o$,\\
                      & $\x(t_f) = \x_f$, \\
                      & $\Delta v(0) = 0$, \\
                      & $d\begin{bmatrix}
                          \x \\ \Delta v
                      \end{bmatrix} = \begin{bmatrix}
                          \mathbf{f}(\x) \\
                          0
                      \end{bmatrix} dt + \begin{bmatrix}
                          \mathbf{0}_3 \\ \hat{\uc} \\ 1
                      \end{bmatrix} d\Gamma$.
    \end{tabular}
    \end{cases}
 \end{align}

Two approaches for Problem $P_\Gamma$ in Equation~\eqref{eq:measure_opt_prob} will be used simultaneously. The first is a piecewise approach, where the measure $\Gamma$ is treated as a sum of impulses and the problem is split into several differentiable segments linked by impulsive discontinuities. The second is the indirect approach, the application of an extended version of Pontryagin's Maximum Principle to impulsive systems~\cite{arutyunov_optimal_2019}. Both will be exposed and applied, as their usage is complimentary.

\section{Piecewise treatment of impulsive maneuvers}\label{sec:piecewise}

In this Section, the bounded variation state trajectory introduced in Section~\ref{sec:orb_man} will be split into a piecewise smooth trajectory, with \emph{impulses} being treated as discontinuities, and the smooth, natural trajectories between impulses being named \emph{coasting arcs}. They can be formally defined as follows.

\begin{definition}[Impulses and Coasting Arcs]
    Let $\texttt{C}$ denote a largest possible interval $[t_c, t_{c+1}]$, $t_c,\ t_{c+1} \in [0, t_f]$, such that
    \begin{equation}
        \Gamma([t_c^+, t_{c+1}^-]) = 0,
    \end{equation}
    and name it a \emph{coasting arc}. Let also $\texttt{I}$ denote a discontinuity in velocity at time $t_i \in [0, t_f]$ at the boundary of a coasting arc interval given by
    \begin{equation}
        \vel(t_i^+) - \vel(t_i^-) = \Gamma(\{t_i\}) \hat{\uc}(t_i),
    \end{equation}
    and name it an \emph{impulse}.
\end{definition}





In other words, solutions to Problem $P_\Gamma$ can be viewed as an alternating sequence of impulses and coasting arcs. The $i$-th impulse is described by the velocity change magnitude \(\Delta v_i = \lVert \vel(t_i^+) - \vel(t_i^-) \rVert \geq 0\) and direction \(\hat{\uc}_i = \hat{\uc}(t_i)\in S^2\). The $c$-th coasting arc is described by its time duration \(d_c = t_{c+1} - t_c > 0\), during which the uncontrolled orbit is propagated.

In the following we consider alternating sequences of impulses and coasting arcs such that the starting time, and end time of each coasting arc is either $0$, $t_f$ or an impulse time $t_i$.
A notation for types of maneuvers described in this piecewise fashion is introduced. \texttt{C} shall represent a coasting arc, and \texttt{I}, an impulse. Any alternating sequence \(\mathcal{S} \in \{\texttt{C}, \texttt{I}\}^{l}\), $l \geq 3$, with 2 or more impulses makes a valid maneuver type: \texttt{ICI}, \texttt{ICIC}, \(\texttt{CICICICIC} = \texttt{C}(\texttt{IC})^4\), etc. Also, let \(\mathcal{I} = \{i \in 1,\dots,l | \mathcal{S}_i = \texttt{I}\}\) and \(\mathcal{C} = \{c \in 1,\dots,l | \mathcal{S}_c = \texttt{C}\}\) be the sets of impulse and coast indices respectively and let $n_i = | \mathcal{I}|$ and $n_c = |\mathcal{C}|$. Sets of feasible maneuvers can be defined as follows.

\begin{definition}
Denote \(M(\mathcal{S}, \x_o, \x_f, t_f)\), or \(M(\mathcal{S})\) for brevity, the set of maneuvers that take the satellite from \(\x_o\) to $\x_f$  in time $t_f$ with a maneuver sequence $\mathcal{S}$. Let $M = \bigcup_{l \geq 3} \bigcup_{\mathcal{S} \in \{\texttt{C}, \texttt{I}\}^{l}} M(\mathcal{S})$ be the set of all possible maneuvers between these endpoints.
\end{definition}

Some definitions relating to the optima of maneuvering problems will also be useful in what follows.

\begin{definition}
Let $\Delta v^\star(\mathcal{S})$ be the optimal cost of maneuvers in the set $M(\mathcal{S})$, and the maneuver that achieves it is called \(M(\mathcal{S})\)-optimal. Also, let $\mathcal{S}^\star$ be the optimal sequence that attains the minimum cost, i.e. $\Delta v^*(\mathcal{S}^\star) = \Delta v^\star$, and the maneuver that achieves it is called $M$-optimal.
\end{definition}

Adding an impulse-coasting arc pair configures a relaxation of the maneuver set. The maneuver set with one less impulse-coasting arc can be recovered by setting the impulse magnitude and the coast duration to zero. Suppose a sequence $\mathcal{S}$ such that $\mathcal{S}_l = \texttt{C}$. Then, it follows immediately from the relaxation that

\begin{align}\label{eq:set_hierarchy}
    M&(\texttt{ICI}) \subset M(\mathcal{S}) \subset M(\mathcal{S}\texttt{IC}) \subset M, \\
    \Delta v^\star&(\texttt{ICI}) \geq \Delta v^\star(\mathcal{S}) \geq \Delta v^\star(\mathcal{S}\texttt{IC}) \geq \Delta v^\star,\label{eq:cost_hierarchy}
\end{align}

The case $\mathcal{S}_l = \texttt{I}$ is analogous. Sequences with more impulses or coasts than $\mathcal{S}^\star$ will have degenerate elements with zero magnitude or duration and have the same cost. So, a method to determine whether a strict cost inequality is found when adding an impulse is required. The primer vector method is, so far, the best fit for this demand. It can indicate whether a given maneuver is a local optimum in a larger maneuver space. For example, a given \texttt{ICI} maneuver can be a local optimum in the $M(\texttt{CICIC})$ set, as opposed to the case where adding infinitesimal coasts at the start and end of the maneuver would lower the cost. It is worth noting that \texttt{ICI} maneuvers correspond to Lambert problems, which are widely studied~\cite{battin_elegant_1984, arora_fast_2014, arora_partial_2015, russell_solution_2019} and often used in orbital maneuvering~\cite{ottesen_unconstrained_2021, sarli_orbit_2015}.

\section{Impulsive optimality condition: primer vector}\label{sec:pv_intro}

Here, the classical primer vector theory will be re-derived with the more appropriate control measure formalism, and different calculation methods for primer vector trajectories will be proposed depending on the dynamical system being studied, with a novel contribution being made in the case of a non-conservative orbital model. 

The term "primer vector" was coined as an analogy with the fact that it imposes a necessary condition for firing the engines, thus acting as a "primer"~\cite{lawden_optimal_1963}. This theory is explored for finite thrust problems in ~\cite{conway_spacecraft_2010}. Typically, the impulsive case is treated as a limiting case where the maximum allowed thrust goes to infinity, and the thrust duration is said to, without proof, reduce to zero. The analysis begins with the impulsive Maximum Principle, stated in Theorem~\ref{theo:optimality_condition}, which is the application of Theorem 2.1 in~\cite{arutyunov_optimal_2019} to orbital maneuvering.

\begin{theorem}~\label{theo:optimality_condition}
    Let $\Gamma^\star, \hat{\uc}^\star, \x^\star, \Delta v^\star$ be a locally optimal solution to Problem $P_\Gamma$ (Equation~\eqref{eq:measure_opt_prob}. Then, there exists an absolutely continuous costate trajectory \(\bPsi: [0, t_f] \rightarrow \R^7\) with components \(\bPsi = \begin{bmatrix}
    \bPsi_r^T & \bPsi_v^T & \Psi_{\Delta v}
\end{bmatrix}^T\), $\bPsi$ not identically null, with which the Hamiltonian function
\begin{equation}
    H(\x, \hat{\uc}, \bPsi) = \bPsi_r^T \vel + \bPsi_v^T \acc(\pos, \vel),
\end{equation}
and the analogous function for the impulsive control
\begin{equation}
    Q(\hat{\uc}, \bPsi) = \bPsi_v^T \hat{\uc} + \Psi_{\Delta v}
\end{equation}
are defined, such that
\begin{align}
    \dot{\bPsi}_s(t) &= - \partial_s H(\x^\star(t), \hat{\uc}^\star(t), \bPsi(t)), s \in \{\x, \Delta v\},\label{eq:psi_dot} \\
    \sup_{\hat{\uc} \in S^2} Q(\hat{\uc}, \bPsi(t)) &\leq 0,\ \forall t \in [0, t_f],\label{eq:supQ} \\
    Q(\hat{\uc}^\star(t), \bPsi(t)) &= 0,\ \Gamma^\star\text{-almost everywhere}.\label{eq:Q0_gae}
\end{align}
\end{theorem}

This theorem follows naturally from the more general Maximum Principle stated with the non-smooth analysis formalism for optimal impulsive control problems ~\cite{arutyunov_optimal_2019}.

These optimality conditions are now analyzed. Remark that
\begin{equation}
    \dot{\Psi}_{\Delta v} = 0,
\end{equation}
and set $\Psi_{\Delta v}(t) = \Psi_{\Delta v} \in \R,\ \forall t \in [0, t_f]$. In Equation~\eqref{eq:supQ}, the supremum is achieved with
\begin{equation}
    \hat{\uc} = \frac{\bPsi_v(t)}{\lVert \bPsi_v(t) \rVert},
\end{equation}
if $\bPsi_v \neq \mathbf{0}_3$, imposing that
\begin{equation}
    \lVert \bPsi_v(t) \rVert \leq - \Psi_{\Delta v},
\end{equation}
at all $t \in [0, t_f]$. From Equation~\eqref{eq:Q0_gae}, \(Q(\hat{\uc}^\star, \bPsi) = 0\) \(\Gamma^\star\)-almost everywhere. Therefore, the support of $\Gamma^\star$ corresponds to instants where \(Q\) is maximum, relating $\bPsi_v$ and $\Psi_{\Delta v}$ by
\begin{equation}\label{eq:psiv_gamma_support}
    \lVert \bPsi_v(t) \rVert = - \Psi_{\Delta v}, \;\Gamma^\star\text{-almost everywhere}.
\end{equation}
 
Introduce now the smooth \textit{primer vector} trajectory $\pv: \R \rightarrow \R^3$ given by
\begin{equation}
    \pv(t) = -\frac{\bPsi_v(t)}{\Psi_{\Delta v}},\ \forall t \in [0, t_f].
\end{equation}

The norm and direction of \(\pv\) indicate when to fire the engines and the optimal impulse direction \(\hat{\uc}\). From now on, it will be assumed that the support of $\Gamma^\star$ is discrete, based on the argument ahead. 
\begin{theorem}\label{theo:isol_zeros}
    Let $\acc(\pos, \vel)$ be an analytic vector field. Then, the solution $\bPsi_x(t) = \begin{bmatrix}
        \bPsi_r^T(t) & \bPsi_v^T(t)
    \end{bmatrix}^T$ to the first six components of Equation~\eqref{eq:psi_dot}, given by
    \begin{equation}
        \begin{bmatrix}
            \dot{\bPsi}_r \\ \dot{\bPsi}_v
        \end{bmatrix} = \begin{bmatrix}
        \mathbf{0}_{3\times 3} & -\left[\frac{\partial \acc(\pos, \vel)}{\partial \pos}\right]^T \\ -\mathbf{I}_3
         & -\left[\frac{\partial \acc(\pos, \vel)}{\partial \vel}\right]^T
    \end{bmatrix} \begin{bmatrix}
        \bPsi_r \\ \bPsi_v
    \end{bmatrix},
    \end{equation}
    is also analytic. Then, any non-constant analytic function of $\bPsi_x$ has isolated zeros.
\end{theorem}
Note that $\bPsi_v^T(t) \bPsi_v(t) - \Psi_{\Delta v}^2$ is analytic and, if not constant, has isolated zeros as stated in Theorem~\ref{theo:isol_zeros}. Therefore from Equation~\eqref{eq:psiv_gamma_support} the support of $\Gamma^\star$ is either the entire interval $[0, t_f]$, or discrete points along it. The first case will be called measure-singular and henceforth ignored. The case where $\Psi_{\Delta v} = 0$ will be called control-singular and similarly ignored, since it is not commonly found in practice~\cite{morelli_characterization_2023}. The necessary conditions for an optimal trajectory from Theorem~\ref{theo:optimality_condition} are restated based on the primer vector in Theorem~\ref{theo:pv_opt}.
\begin{theorem}~\label{theo:pv_opt}
    If a solution to Problem $P_\Gamma$ is locally optimal, then its primer vector trajectory $\pv$ satisfies
\begin{enumerate}
    \item \(\mathbf{p}(t)\) and \(\dot{\mathbf{p}}(t)\) are continuous for all \(t \in [0, t_f]\);
    \item \(\lVert \mathbf{p}(t) \rVert \leq 1\) for all \(t \in [0, t_f]\);
    \item \(\mathbf{p}(t) = \hat{\uc}(t)\) at the impulse instants, with the corollary that \(\lVert \mathbf{p}(t) \rVert = 1\);
    \item \(\frac{d \lVert \mathbf{p} \rVert}{dt} = 0\) at impulses at times \(t \in (0, t_f)\).
\end{enumerate}
\end{theorem}

These conditions are widely found in the literature~\cite{conway_spacecraft_2010,luo_interactive_2010, lion_primer_1968, jezewski_efficient_1968} for the finite thrust case. The impulsive thrust case is usually treated as a limiting case, supported by the intuition that burn times become smaller when the maximum thrust grows. However, they had never been derived with the more correct impulsive optimal control formalism, which justifies the presentation of this analysis.

\subsection{Primer vector algorithm}\label{sec:pv_alg}

This Section shall explore how the primer vector may be used to improve trajectories, as well as how to compute it for different classes of orbital models.

Theorem~\ref{theo:optimality_condition} gives optimality conditions for the \(M\)-optimal maneuver. \(M(\mathcal{S)}\)-optimal maneuvers, $\mathcal{S} \neq \mathcal{S}^\star$, may violate them. Lion and handelsman~\cite{lion_primer_1968} expanded on Lawden's theory, developing  modifications for the sequence \(\mathcal{S}\) based on the primer vector trajectory of the \(M(\mathcal{S)}\)-optimal maneuver from a first order calculus of variations analysis:
\begin{algorithm}
\caption{Primer vector algorithm.}
\begin{algorithmic}
    \STATE $\mathcal{S} \gets \texttt{ICI}$
    \STATE $\pv \gets $ primer vector trajectory of $M(\mathcal{S})$-optimal trajectory
    \IF{$\mathcal{S}_1 = \texttt{I}$ \AND \(\lVert \mathbf{p}(0) \rVert = 1\) \AND \(\partial_t \lVert \mathbf{p} \rVert \mid_{t=0} > 0\)}
    \STATE $\mathcal{S} \gets \texttt{C}\mathcal{S}$ \COMMENT{Add initial coast}
    \ENDIF
    \IF{$\mathcal{S}_l = \texttt{I}$ \AND \(\lVert \mathbf{p}(t_f) \rVert = 1\) \AND \(\partial_t \lVert \mathbf{p} \rVert \mid_{t=t_f} < 0\)}
    \STATE $\mathcal{S} \gets \mathcal{S}\texttt{C}$ \COMMENT{Add final coast}
    \ENDIF
    \STATE $\pv \gets $ primer vector trajectory of $M(\mathcal{S})$-optimal trajectory
    \WHILE{\(\lVert \mathbf{p}(t) \rVert > 1\) for some \(t \in [0, t_f]\)}
    \STATE $\mathcal{S} \gets \texttt{CI}\mathcal{S}$ \COMMENT{Add mid-course impulse}
    \STATE $\pv \gets $ primer vector trajectory of $M(\mathcal{S})$-optimal trajectory
    \ENDWHILE
    \STATE $\mathcal{S}^\star \gets \mathcal{S}$
\end{algorithmic}
\end{algorithm}

Iterating over these rules to search \(\mathcal{S}^\star\) has often been used in the literature for the Keplerian model ~\cite{luo_interactive_2010, conway_spacecraft_2010, prussing_optimal_1986, taheri_how_2020, lion_primer_1968, jezewski_efficient_1968}. Now, some means of computing the primer vector is searched. 

Computing the primer vector trajectory is done piecewise on coasting arcs between two consecutive impulses, since the primer vector is known at impulse times. This configures a linear two point boundary value problem (LTPBVP) ~\cite{conway_spacecraft_2010}. Should there be a coasting arc not bounded by two impulses as in \texttt{C}$\dots$ or $\dots$\texttt{C}, the neighboring primer vector solution is propagated. Suppose two consecutive impulses happen at times $t_1$, $t_2$ such that $0 < t_1 < t_2 < t_f$. The primer vector values at these instants are given by:
\begin{align}\label{eq:pv_deltav}
    \mathbf{p}(t_1) &= \hat{\uc}(t_1),\\ 
    \mathbf{p}(t_2) &= \hat{\uc}(t_2).
\end{align}

This LTPBVP will now be approached with a general non-conservative model first. The optimality condition stated in Equation~\eqref{eq:psi_dot} gives rise to second order linear primer vector dynamics with dynamic matrix \(\mathbf{A}_p(\pos, \vel)\):
\begin{equation}\label{eq:pdot_ncon}
    \begin{bmatrix}
        \dot{\mathbf{p}} \\ \ddot{\mathbf{p}}
    \end{bmatrix} = \begin{bmatrix}
        \mathbf{0}_{3 \times 3} & \mathbf{I}_3 \\
        \left[\frac{\partial \acc(\pos, \vel)}{\partial \pos}\right]^T & -\left[\frac{\partial \acc(\pos, \vel)}{\partial \vel}\right]^T
    \end{bmatrix} \begin{bmatrix}
        \mathbf{p} \\ \dot{\mathbf{p}}
    \end{bmatrix} = \mathbf{A}_p(\pos, \vel) \begin{bmatrix}
        \mathbf{p} \\ \dot{\mathbf{p}}
    \end{bmatrix}.
\end{equation}

Therefore, the primer vector admits a primer vector transition matrix (PVTM) \(\bPhi_p(t)\):
\begin{equation}\label{eq:pvtm}
    \begin{bmatrix}
        \mathbf{p}(t) \\ \dot{\mathbf{p}}(t)
    \end{bmatrix} = \mathbf{\mathbf{\Phi}}_p(t - t_1) \begin{bmatrix}
        \mathbf{p}(t_1) \\ \dot{\mathbf{p}}(t_1)
    \end{bmatrix}.
\end{equation}

If \(\bPhi_p(t - t_1)\) is known, Equation~\eqref{eq:pvtm} applied to \(t=t_2\) allows for the determination of \(\dot{\pv}(t_1)\) and consequently, of the entire primer vector trajectory. The most straightforward method of calculating it, albeit not present in the primer vector literature, is simply using the initial value problem in Equation~\eqref{eq:stm_ivp}, which is to be integrated alongside the the equations of motion, which will be named the ODE method of PVTM calculation:
\begin{equation}\label{eq:PVTM_ODE}
    \begin{cases}
        \dot{\mathbf{\mathbf{\Phi}}}_p(t - t_1) = \mathbf{A}_p(\pos(t), \vel(t)) \mathbf{\mathbf{\Phi}}_p(t-t_1), t\in[t_1, t_2] \\
        \mathbf{\mathbf{\Phi}}_p(0) = \mathbf{I}_6
    \end{cases}.
\end{equation}

If the model being treated is conservative, \(\partial_{\vel} \acc = \mathbf{0}_{3 \times 3}\) and \(\partial_\pos \acc = [\partial_\pos \acc]^T\) since \(\partial_\pos \acc\) is actually the Hessian of the potential field. In this case, \(\mathbf{A}_p(\pos) = \partial_\x \mathbf{f}(\x)\), the Jacobian of the dynamics, and the PVTM and the STM are equal:
\begin{equation}
    \bPhi_p(t - t_1) = \bPhi_\delta(t - t_1) = \frac{\partial \x(t)}{\partial \x(t_1)},
\end{equation}
which is convenient if a practical mean of calculation the sensitivity of the orbital trajectory is available. This will be named the STM method of PVTM calculation.

If in addition the model is Keplerian, the analytical expression for the STM provided by Glandorf~\cite{glandorf_lagrange_1969} can be used, which will be named the Glandorf method of PVTM calculation. From Equation~\eqref{eq:stm_glandorf},
\begin{equation}
    \bPhi_p(t - t_1) = \bPhi_\delta(t - t_1) = \mathbf{P}(t) \mathbf{P}^{-1}(t_1).
\end{equation}

The Glandorf method is commonly found in the literature, as in \cite{conway_spacecraft_2010, prussing_optimal_1986}, with a similar treatment found in~\cite{lion_primer_1968}, and some works use the primer vector through the STM method in the Keplerian or other conservative models, including multi-body problems~\cite{bokelmann_optimization_2020, jezewski_efficient_1968, bucchioni_optimal_2023}. However, the application of primer vector to non-conservative models with the ODE method is, to the authors knowledge, novel. Table~\ref{tab:pv_calc} summarizes the domain of validity of each method according to the class of orbital model used. 

\begin{table}[htbp]
    \centering
    \begin{tabular}{cccc} \toprule
        Model Class & Glandorf & STM & ODE \\ \midrule
        Keplerian & \checkmark & \checkmark & \checkmark \\
        Conservative & \(\times\) & \checkmark & \checkmark \\
        Non-conservative & \(\times\) & \(\times\) & \checkmark \\ \bottomrule
    \end{tabular}
    \caption{Summary of primer vector calculation methods for different model types.}
    \label{tab:pv_calc}
\end{table}

\section{Optimization framework}\label{sec:impl}

The numerical implementation of this work shall be discussed in the present Section. The piecewise approach from Section~\ref{sec:piecewise} is used as the basis for a multiple shooting scheme. An eighth-order Runge-Kutta (RK8) method~\cite{verner_numerically_2010} was used for discretizing the dynamical equations. Let $\mathbf{f}_{RK}$ be the RK8 integrator such that
\begin{equation}
    \x(t + \Delta t) = \mathbf{f}_{RK}(\Delta t, \x(t))
\end{equation}
for some time step $\Delta t$ and discretization instant $t$.

The input parameters are:
\begin{enumerate}
    \item \(\mathbf{r}_1\), \(\mathbf{v}_1\): initial orbital position and velocity;
    \item \(\mathbf{r}_2\), \(\mathbf{v}_2\): final orbital position and velocity;
    \item \(t_f\): transfer time;
    \item \(N\): number of integration steps per coasting arc.
    \item \(\mathcal{S} \in \{\texttt{C}, \texttt{I}\}^{n_c + n_i}\): the desired sequence of coasts and impulses.
\end{enumerate}

The problem's definition and variable set dynamically depends on the input sequence. The \(i\)-th impulse is described by the variables:
\begin{enumerate}
    \item \(\Delta v_i \geq 0\): magnitude of the impulse;
    \item \(\hat{\uc}_i \in \R^3\): direction of the impulse,
\end{enumerate}
which is a parameterization of impulses that avoids norm singularities, and follows naturally from Problem $P_\Gamma$ in Equation~\eqref{eq:measure_opt_prob}.

The \(c\)-th coasting arc is described by:
\begin{enumerate}
    \item \(d_c \geq 0\): total duration of the arc;
    \item \(\x^j_c \in \R^6, j=1,\dots,N\): state vector variables for each coasting arc. 
\end{enumerate}

 The full multiple shooting problem for the optimization of the impulsive maneuver with sequence \(\mathcal{S}\), with $N$ integration steps, Problem $P^N_\mathcal{S}$, is:

\begin{align}\label{eq:ms_problem}
    \begin{tabular}{cl}
     \(\Delta v^\star(\mathcal{S}) = \min\)                              & \(\sum_{i \in \mathcal{I}} \Delta v_i\)\\
    \(\Delta v_i  \in \R_+, i \in \mathcal{I}\), &  \\
    \(\hat{\uc}_i \in \R^3, i \in \mathcal{I}\),     & \\
    \(d_c         \in \R_+, c \in \mathcal{C}\),     & \\
    \(\x^j_c      \in \R^6, j=1,\dots,N, c \in \mathcal{C}\) & \\
    & \\
    \textbf{subject to:}        & \\
    & \\
    Total time                  & \(\sum_{c \in \mathcal{C}} d_c = t_f\) \\
    Unit directions             & \(\hat{\uc}_i^T \hat{\uc}_i = 1, i \in \mathcal{I}\) \\
    Propagation of coasts       & \(\x^{j+1}_c = f_{RK}(\frac{d_c}{N-1}, \x^j_c)\), \\
                                & \(j=1,\dots,N-1, c \in \mathcal{C}\) \\
    Impulse boundary conditions & \(\x^1_{i+1} = \x^{N}_{i-1} + \begin{bmatrix}
        0_{3\times1} \\ \Delta v_i \hat{\mathbf{u}}_i
    \end{bmatrix}\), \\
                                & \(i \in \mathcal{I}/\{1, n_i+n_c\}\) \\
    Initial condition           & \(\begin{cases}
                                        \x_1^1 = \begin{bmatrix}
                                            \mathbf{r}_1 \\ \mathbf{v}_1
                                        \end{bmatrix}, 1 \in \mathcal{C} \\
                                        \x_2^1 = \begin{bmatrix}
                                            \mathbf{r}_1 \\ \mathbf{v}_1 + \Delta v_1 \hat{\mathbf{u}}_1
                                        \end{bmatrix}, 1 \in \mathcal{I}
                                    \end{cases}\) \\
    Final condition             & \(\begin{cases}
                                        x_c^N = \begin{bmatrix}
                                            \mathbf{r}_2 \\ \mathbf{v}_2
                                        \end{bmatrix}, c = n_i+n_c \in \mathcal{C} \\
                                        x_{i-1}^N + \begin{bmatrix}
                                            \mathbf{0}_{3} \\ \Delta v_i \hat{\mathbf{u}}_i
                                        \end{bmatrix} = \begin{bmatrix}
                                            \mathbf{r}_2 \\ \mathbf{v}_2
                                        \end{bmatrix}, i = n_i+n_c \in \mathcal{I}
                                    \end{cases}\) \\
                                & \(\)
    \end{tabular}
 \end{align}

The numerical approach to this work was developed in Julia~\cite{bezanson_julia_2017}, with the maneuver optimization done in a multiple shooting scheme implemented in CasADi~\cite{andersson_casadi_2019} using the nonlinear solver Ipopt~\cite{wachter_implementation_2006}. The orbital dynamics functions were implemented to be compatible with the symbolic types introduced by CasADi, and the implementation was validated with the SatelliteToolbox package~\cite{chagas_juliaspacesatellitetoolboxjl_2025}, developed at Brazil’s National Institute of Space Research (INPE), which does not directly participate in the optimization.

The problem of optimizing orbital maneuvers is non-convex, and has many local optima, possibly infinitely many~\cite{saloglu_existence_2023}. Ensuring the quality of solutions to nonconvex problems is notoriously hard. Evolutionary algorithms are an increasingly promising solution to this~\cite{luo_interactive_2010}, but they require large amounts of function evaluations even for simple cases. Multiple shooting schemes are amenable to gradient-based algorithms, and can be solved with modern tools such as SNOPT~\cite{gill_snopt_2005} or Ipopt, which is not uncommon in the field of orbital maneuvering~\cite{bokelmann_optimization_2020}~\cite{russell_primer_2007}.  

To offset the local nature of this method, a number of random initial guesses can be generated. For scenarios where the starting and final orbits are similar, setting impulse magnitudes to zero has been found to lead to workable guesses. The most important degrees of freedom in the initial guess choice are coast durations, which can be uniformly sampled as detailed in Appendix~\ref{app:unif_sampl}. 

Many methods for integrating dynamical systems for optimal control exist, such as collocation, pseudospectral methods, single shooting, and multiple shooting~\cite{betts_practical_2010}. Multiple shooting was chosen since it leads to a numerically stable problem and is standard in optimal control. This discretization scheme is agnostic to the dynamics of the system, that is, one can "drop-in" the Keplerian model, the J2 model, or the J2-Drag model and the scheme remains the same, which is quite convenient.

The issues of variable normalization and a smooth approximation of the atmospheric density model are discussed in Appendix~\ref{app:var_norm} and~\ref{app:smooth_atm} respectively. Problem $P^N_\mathcal{S}$ in Equation~\eqref{eq:ms_problem} is repeatedly solved for sequences $\mathcal{S}$ with increasing number of impulses by iterating over the rules described in Section~\ref{sec:pv_alg}.

\section{Numerical Results}\label{sec:res}

This Section presents the scenarios of orbital maneuvering chosen to illustrate the primer vector method before their solutions under different orbital models are given. All necessary orbital parameters (\(\mu\), \(J_2\), \(R\), \(\omega_E\)) were set to the values offered by the \texttt{SatelliteToolbox}. Two scenarios were chosen, one where a coplanar transfer between two circular orbits is done in half a revolution, and a noncoplanar rendez-vous case spanning two full revolutions. Both cases are presented in detail ahead, and summarized in Table~\ref{tab:scenario_orb_elems}.

In the first scenario, the transfer time was set so that the Hohmann transfer is recovered in the Keplerian model. This is a simple scenario that highlights the difference between the different orbital models used in the work, and validates the developed optimization setup. A second, more challenging scenario is taken from \cite{luo_interactive_2010}, where a noncoplanar rendez-vous scenario is proposed and solved under the Keplerian model. The transfer time is equal to two orbital periods of the target orbit. Again, the goal is validating the code against a published result under the Keplerian model, and exploring how the optimal trajectory may be different under other more realistic models. 

\begin{table}[htpb]
    \centering
        \begin{tabular}{lcccc} \toprule
                   & \multicolumn{2}{c}{Circle to Circle} & \multicolumn{2}{c}{Noncoplanar rendez-vous} \\ \midrule
        Element    & Initial        & Final              & Initial & Final \\ \midrule
        Semi-major axis      & \(7000.0\) km  & \(9000.0\) km      & \(6748.1\) km         & \(6778.1\) km   \\
        Eccentricity      & \(0.0\)        & \(0.0\)            & \(0.0\)            & \(0.0\)        \\
        Inclination      & \(51.0^\circ\) & \(51.0^\circ\)     & \(42.1^\circ\)      & \(42.0^\circ\) \\
        RAAN & \(0.0^\circ\)  & \(0.0^\circ\)      & \(120.2^\circ\)   & \(120.0^\circ\)  \\
        Argument of periapsis & \(0.0^\circ\)  & \(0.0^\circ\)      & \(0.0^\circ\)  & \(0.0^\circ\)  \\
        True anomaly & \(0.0^\circ\)  & \(180.0^\circ\)    & \(175.0^\circ\)  & \(180.0^\circ\)  \\ 
        \textbf{Maneuver time} & \multicolumn{2}{c}{\(3560.541\) s} & \multicolumn{2}{c}{\(11107.158\) s} \\\bottomrule
    \end{tabular}
    \caption{Orbital elements used in each maneuver scenario.}
    \label{tab:scenario_orb_elems}
\end{table}

The primer vector algorithm was applied to the described scenarios. Both of them were solved under the three orbital models summarized in Table~\ref{tab:orb_models}, and numerical results will be presented in this section. The Keplerian model results are left for the summarized discussion at the end, since they correspond to the notable analytical Hohmann transfer and to published results~\cite{luo_interactive_2010} and were used for code validation. Remarkable differences between the Keplerian and J2 maneuvers were found, and the modelling of drag was found to have little impact on trajectories. In all primer vector plots in this section, vertical dashed lines indicate impulse times.

\subsection{J2 model}

Optimized maneuvers are presented for the aforementioned scenarios under the J2 model introduced in Section~\ref{sec:orb_models}. Under this model, and according to the theory developed in Section~\ref{sec:pv_alg}, only the STM and ODE methods correctly calculate primer vector trajectories, and both will be presented for all cases. Early testing showed that they differ from the Glandorf method, not applicable to this model, as expected.

\subsubsection{Circle to Circle Scenario}

\begin{figure}
    \centering
    \includegraphics[width=\linewidth]{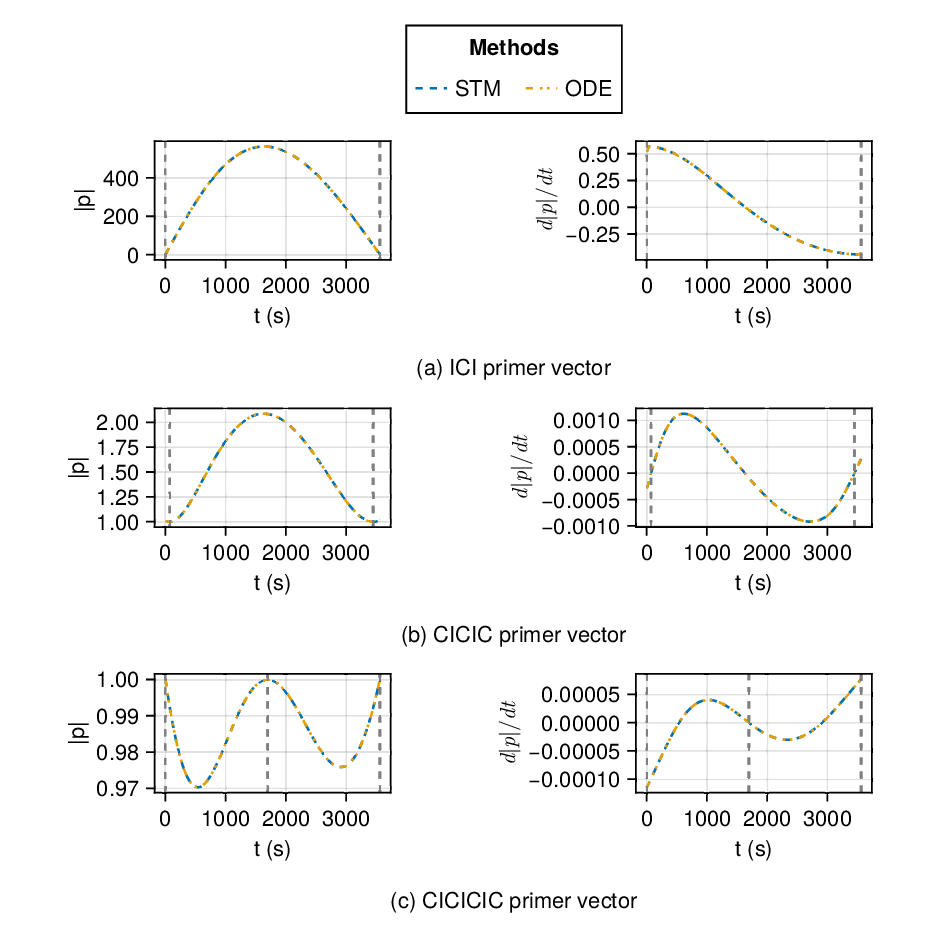}
    \caption{Primer vector histories of circle to circle rendez-vous under the J2 model.}
    \label{fig:j2_c2c}
\end{figure}

The simple circle to circle rendez-vous case was first optimized with a \texttt{ICI} maneuver, followed by iteration over the primer vector rules in Section~\ref{sec:pv_alg}, resulting in a three impulse maneuver. Primer vector trajectories can be seen in Figure~\ref{fig:j2_c2c}, and impulse descriptions, in Table~\ref{tab:j2_c2c}

The \texttt{ICI} solution is very expensive due to large plane changes required for a trajectory connecting antipodal points under the J2 disturbance, which tends to twist trajectories out-of-plane. This solution is not a local optimum in the space of \texttt{CICIC} solutions, as indicated by the non-zero endpoint primer vector norm derivatives in Figure~\ref{fig:j2_c2c} (a), positive at the start and negative at the end, suggesting the addition of initial and final coasts to reduce the very large cost of this solution.

The solution of the \texttt{CICIC} maneuver closely resembles a Hohmann transfer in its impulse magnitudes. But remarkably, the impulses are not parallel to the satellite's velocity, suggesting the need for plane corrections. The primer vector norm exceeds unity, as seen in Figure~\ref{fig:j2_c2c} (b), showing that this trajectory is not a local optimum in the \texttt{CICICIC} maneuver space, so adding another impulse should decrease the cost.

Now, the \texttt{CICICIC} case is solved, with the primer vector trajectory in Figure~\ref{fig:j2_c2c} (c) satisfying optimality conditions. A very small intermediate impulse is made about halfway through the maneuver for plane correction. The necessary conditions on the primer vector trajectory for optimality are satisfied with 3 impulses, and the final cost of the J2 model is about \(6\) m/s higher than the analytical Hohmann transfer solution, for the same maneuver time.

\begin{table}[] \small
    \begin{tabular}{cccccc} \toprule
    \(\mathcal{S}\) &\textbf{Impulse}   & 1 & 2 & 3           & \textbf{Total} (m/s) \\ \midrule
    ICI & \(t\) (s)           & 0.0 & 3560.54 &  &\\
                & \(\Delta v\) (m/s)    & 5209.48 & 4318.72 &  & 9528.2    \\ \midrule
    CICIC & \(t\) (s)           & 72.53 & 3448.64 &  &\\
                & \(\Delta v\) (m/s)    & 473.47 & 438.47 &  & 911.93    \\ \midrule
    CICICIC & \(t\) (s)           & 0.39 & 1697.06 & 3559.91 &\\
                & \(\Delta v\) (m/s)    & 451.26 & 25.17 & 416.62 & 893.05    \\ \bottomrule
    \end{tabular}
    \caption{Impulse descriptions of circle to circle rendez-vous under the J2 model.}\label{tab:j2_c2c}
\end{table}

\subsubsection{Noncoplanar Rendez-Vous Scenario}

\begin{figure}
    \centering
    \includegraphics[width=\linewidth]{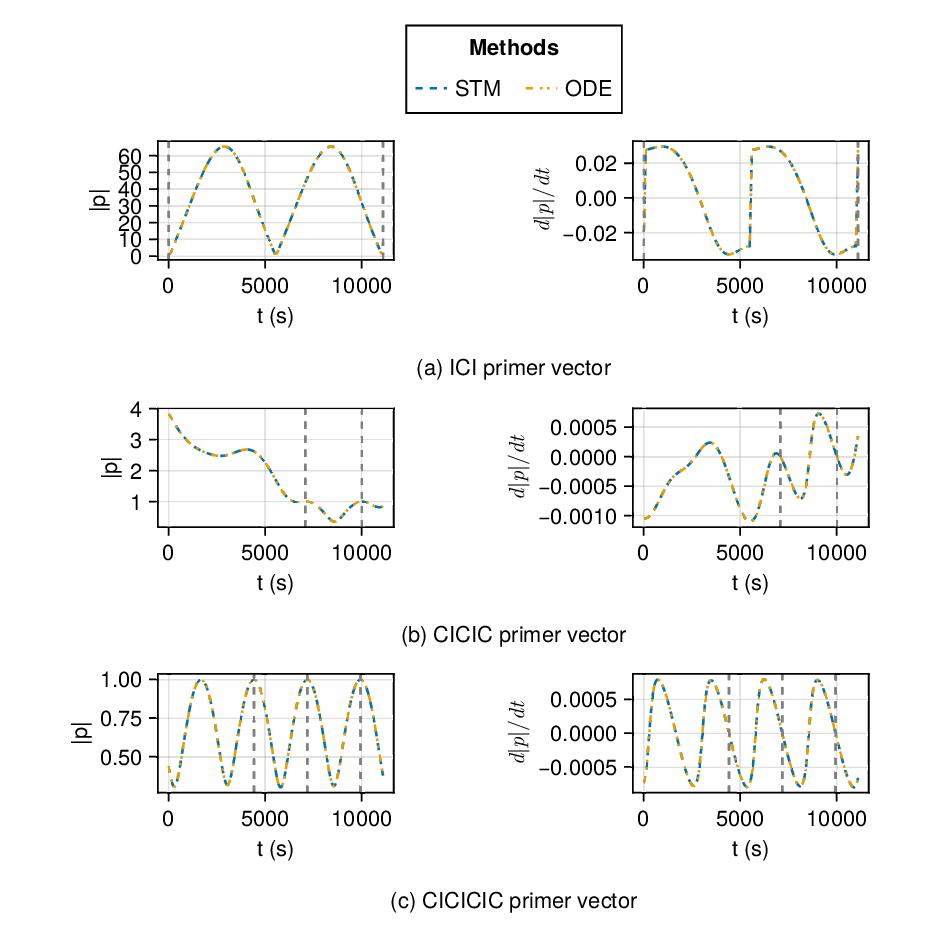}
    \caption{Primer vector histories of noncoplanar rendez-vous under the J2 model.}
    \label{fig:j2_ncop}
\end{figure}

Now, the noncoplanar rendez-vous scenario is tackled under the J2 model, with primer vector trajectories shown in Figure~\ref{fig:j2_ncop} and impulse descriptions in Table~\ref{tab:j2_ncop}. As usual, the first maneuver to be solved is a \texttt{ICI} case. The maneuver is quite expensive, and does not satisfy primer vector necessary conditions. In particular, primer vector norm derivative is not continuous, and the primer vector norm is greater than unity. A \texttt{CICIC} maneuver is searched next.

Next, the \texttt{CICIC} maneuver is only \(5\) m/s more expensive than its Keplerian counterpart (see comparison in Section~\ref{sec:mod_comp}). Despite the small cost difference, impulse times and magnitudes are considerably different to those found in the Keplerian counterpart. The primer vector trajectory is now continuously differentiable, but still violates the unity norm condition. Now, adding an impulse is necessary.

A \texttt{CICICIC} maneuver was solved next. The cost is marginally lower than the previous case, which indicates that only small improvements are available in this case after the 2-impulse solution, but now the primer vector necessary conditions are satisfied. It is known that sometimes following the primer vector algorithm will yield small improvements in cost at the expense of adding very small impulses~\cite{bokelmann_optimization_2020}. It is worth clarifying that the primer vector rules require that \(\lVert \mathbf{p} \rVert = 1\) at every impulse time, but this does not mean that unity norm cannot occur outside impulse times, which exactly what happens at around \(4000\) s in this case.

\begin{table}[] \small
    \begin{tabular}{cccccc} \toprule
    \(\mathcal{S}\) &\textbf{Impulse}   & 1       & 2          & 3          & \textbf{Total} (m/s) \\ \midrule
 \texttt{ICI} & \(t\) (s)          & 0.0          & 11107.1576          & &\\
 & \(\Delta v\) (m/s) & 743.66668& 727.80415& &1471.47082\\ \midrule
 $\texttt{C}(\texttt{IC})^2$ & \(t\) (s)          & 7080.30484          & 10011.10872          & &\\
 & \(\Delta v\) (m/s) & 24.12503& 34.14936& &58.27439\\ \midrule
 $\texttt{C}(\texttt{IC})^3$ & \(t\) (s)          & 1676.61473          & 7185.69293          & 9942.01138          &\\
 & \(\Delta v\) (m/s) & 5.84342& 20.83452& 29.32859&56.00653\\ \bottomrule
    \end{tabular}
    \caption{Impulse descriptions of noncoplanar rendez-vous under the J2 model.}
    \label{tab:j2_ncop}
\end{table}

\subsection{J2+Drag model}

Finally, the non-conservative model of J2 and drag disturbances is tackled. The satellite was supposed to have \(C_D = 2.2\), \(S = \frac{\pi}{4}\ \text{m}^2\), and \(m = 10\) kg, which is representative of LEO Cubesats, largely used nowadays. This model's results will be presented in a summarized fashion, since they were found to be almost identical to J2 results. In addition, although the theory predicts that only the ODE method of primer vector calculation should work for this system, as developed in Section~\ref{sec:pv_alg}, it was found that the STM method gave almost identical trajectories, as will be discussed.

\subsubsection{Circle to Circle Scenario}

\begin{figure}
    \centering
    \includegraphics[width=\linewidth]{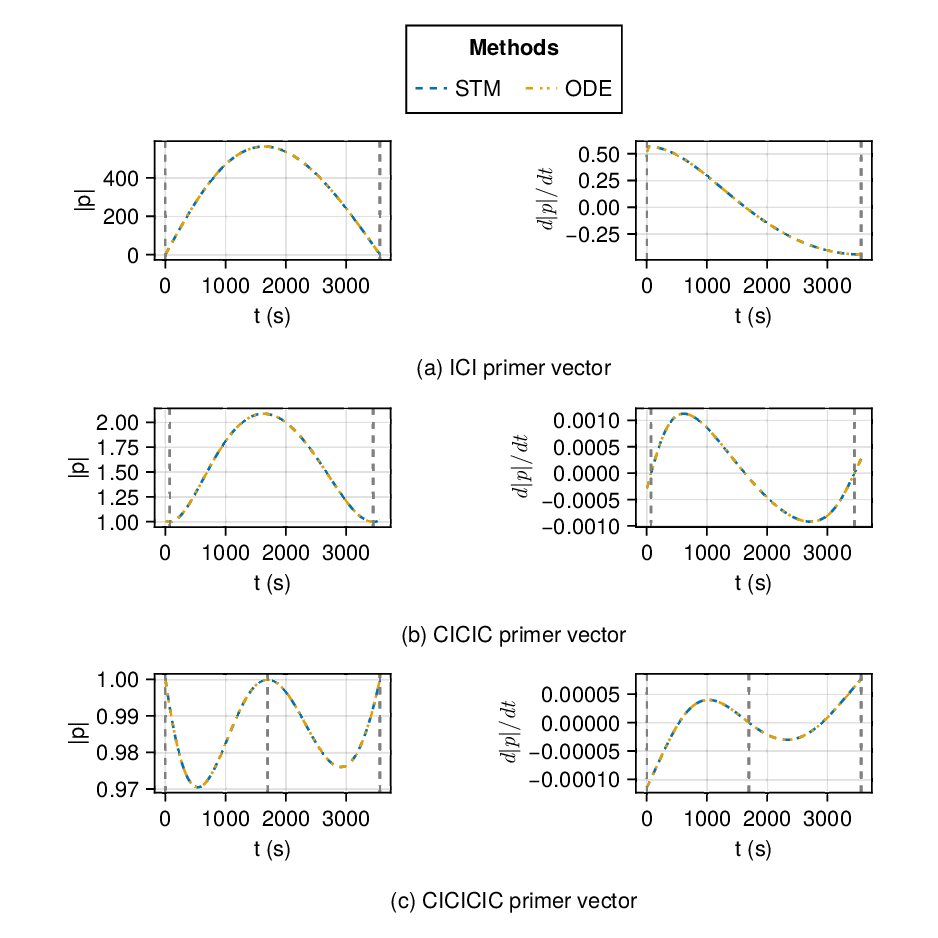}
    \caption{Primer vector histories of circle to circle rendez-vous under the J2+Drag model.}
    \label{fig:j2drag_c2c}
\end{figure}

The circle to circle case was repeatedly solved according to the primer vector algorithm up to three impulses. The maneuvers' primer vector histories and costs are displayed in Figure~\ref{fig:j2drag_c2c} and Table~\ref{tab:j2drag_c2c}. Straightaway, it is clear that these results are almost identical to the ones found in the J2 model, given in Table~\ref{tab:j2_c2c}. This can be easily explained by the fact that the effects of drag are some orders of magnitude weaker than J2 effects at the considered altitudes. 

The theory from Section~\ref{sec:pv_alg} disallows, in general, the use of the STM method with non-conservative models. However, as seen in Figure~\ref{fig:j2drag_c2c}, both STM and ODE methods give almost identical primer vector trajectories. This can be explained by comparing the matrices \(\mathbf{A}_\delta\) and \(\mathbf{A}_p\), the coefficients of the ODEs of the STM and PVTM respectively. Indeed, the matrices satisfy \((\mathbf{A}_\delta)_{\substack{4\leq i \leq 6 \\ 4 \leq j \leq 6}} = - (A_p)_{\substack{4\leq i \leq 6 \\ 4 \leq j \leq 6}}^T\), as predicted by the theory in Section~\ref{sec:pv_alg}, but due to the low intensity of drag forces, they are both close to \(\mathbf{0}_{3\times 3}\), which is the lower right quadrant of \(\mathbf{A}_{p}\) in the J2 model. This shows that for maneuvers where drag is weak, either by virtue of the short timespan or high altitude, the effects of drag may be neglected for primer vector calculation. For the initial state in this scenario, for instance, \(\lVert \mathbf{A}_p - \mathbf{A}_\delta \rVert _\infty < 10^{-11}\), proving that indeed drag can be neglected for maneuvers of this type.

\begin{table}[] \small
    \begin{tabular}{cccccc} \toprule
    \(\mathcal{S}\) &\textbf{Impulse}   & 1 & 2 & 3           & \textbf{Total} (m/s) \\ \midrule
    ICI & \(t\) (s)           & 0.0 & 3560.5408 &  &\\
                & \(\Delta v\) (m/s)    & 5209.4785 & 4318.7183 &  & 9528.1968    \\ \midrule
    CICIC & \(t\) (s)           & 72.5315 & 3448.6386 &  &\\
                & \(\Delta v\) (m/s)    & 473.4671 & 438.4661 &  & 911.9333    \\ \midrule
    CICICIC & \(t\) (s)           & 0.1416 & 1696.7899 & 3560.3114 &\\
                & \(\Delta v\) (m/s)    & 451.2256 & 25.2157 & 416.5868 & 893.0281    \\ \bottomrule
    \end{tabular}
    \caption{Impulse descriptions of circle to circle rendez-vous under the J2+Drag model.}\label{tab:j2drag_c2c}
\end{table}

\subsubsection{Noncoplanar Rendez-Vous Scenario}

The noncoplanar rendez-vous case with drag was the most challenging problem solved in this work. The complexity of the model, with 28 terms to describe atmospheric density, combined with the sixteen nested dynamics function calls required by the RK8 method, meant that the CasADi models were very large, requiring some tens of gigabytes of memory for their solution.

\begin{figure}
    \centering
    \includegraphics[width=\linewidth]{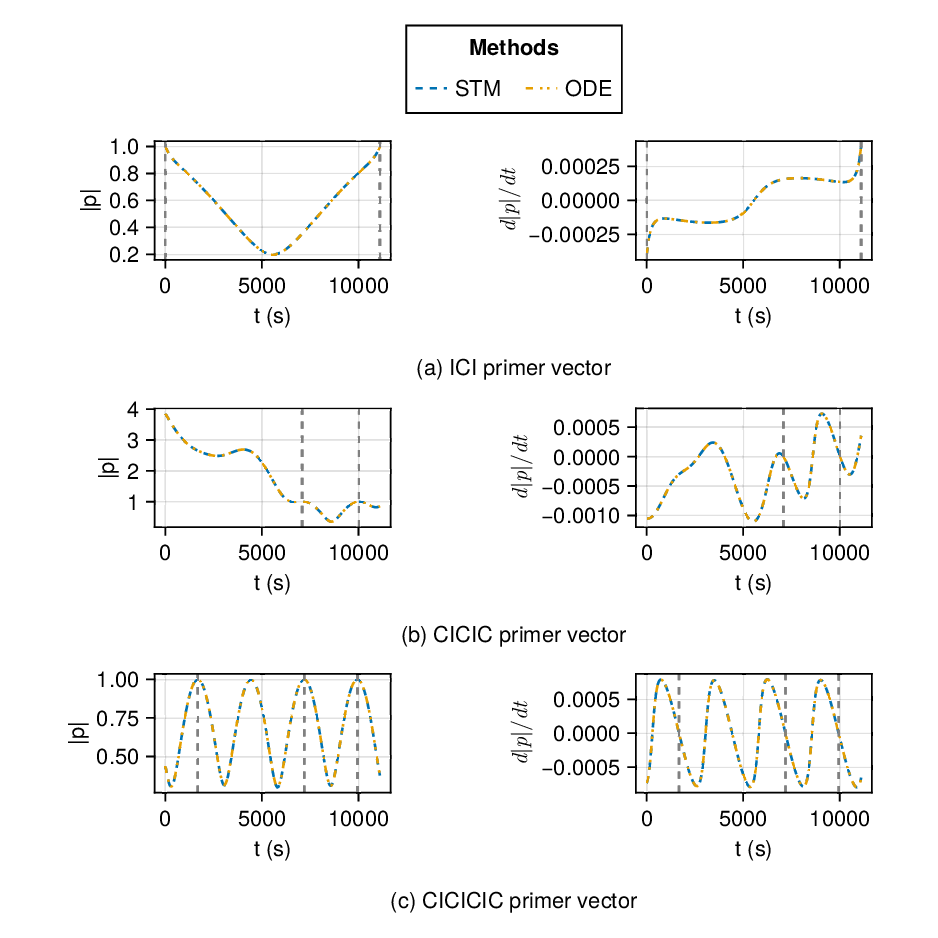}
    \caption{Primer vector histories of noncoplanar rendez-vous under the J2+Drag model.}
    \label{fig:j2drag_ncop}
\end{figure}

The results in this case are presented in summarized fashion in Table~\ref{tab:j2drag_ncop}, along with the primer vector trajectories in Figure~\ref{fig:j2drag_ncop}.

The J2+Drag \texttt{ICI} maneuver has very high cost. It is in a way a ``brute force'' maneuver, where a very eccentric, high semi major axis arc is followed instead of a multi-revolution trajectory. A multi-revolution \texttt{ICI} maneuver, as was found in the J2 model, was searched by tuning the initial guesses, but only solver errors and infeasible results were found, suggesting the conclusion that drag renders that multi-revolution J2 maneuver infeasible. The next two maneuvers were found to be almost identical to the J2 maneuvers, as will be discussed in Section~\ref{sec:mod_comp}.

All J2+Drag primer vector trajectories calculated with the STM and ODE methods are visually identical, and indistinguishable to J2 primer vector trajectories. This shows that the scenarios considered are cases where drag is too small to distinguish the STM and ODE methods, and to greatly distinguish J2+Drag maneuvers from J2 maneuvers. Therefore, it is reasonable to conjecture that both should happen together: a scenario where J2 and J2+Drag maneuvers are sufficiently distinct will also have STM and ODE primer vector trajectories that are sufficiently distinct. A maneuvering scenario including aerobraking or very long time periods could be investigated in future works to validate or refute this conjecture.

\begin{table}[] \small
    \begin{tabular}{cccccc} \toprule
    \(\mathcal{S}\) &\textbf{Impulse}   & 1 & 2 & 3           & \textbf{Total} (m/s) \\ \midrule
    ICI & \(t\) (s)           & 0.0 & 11107.1576 &  &\\
                & \(\Delta v\) (m/s)    & 11743.61 & 11711.0072 &  & 23454.6171    \\ \midrule
    CICIC & \(t\) (s)           & 7078.6538 & 10011.9002 &  &\\
                & \(\Delta v\) (m/s)    & 24.0895 & 34.2266 &  & 58.3161    \\ \midrule
    CICICIC & \(t\) (s)           & 1676.6212 & 7185.7104 & 9941.9901 &\\
                & \(\Delta v\) (m/s)    & 5.8844 & 20.8038 & 29.326 & 56.0142    \\ \bottomrule
    \end{tabular}
    \caption{Impulse descriptions of noncoplanar rendez-vous under the J2+Drag model.}\label{tab:j2drag_ncop}
\end{table}

\subsection{Summary of results}\label{sec:mod_comp}

The discussion is finished with a comparative analysis of each scenario under the different orbital models. Table~\ref{tab:model_comparison} compares the optimal number of impulses and optimal cost for each combination of scenario and model. Here, it is worth highlighting that the circle-to-circle scenario under the Keplerian model corresponds to the well-known Hohmann Transfer.

Despite the similar final costs between the three models for the circle to circle scenario, equivalent maneuver sequences give very different optima, as highlighted by the \texttt{ICI} sequence, which is the Hohmann transfer under the Keplerian model, costs $\Delta v^\star = 9528.19$ m/s under the J2 and J2+Drag models. It is interesting that similar costs can only be found with different number of impulses.

Both scenarios had higher costs under the perturbed models when compared with the Keplerian model. This should not generalize, and is in fact an artifact of the scenario choice. The important takeaway is that costs, instants and number of impulses can be different between maneuvers in the Keplerian and J2 models, even at time frames of 2 revolutions or less.

The very similar costs found in the J2 and J2+Drag models were consistently reproduced with multiple code executions. This means the small gap can truly be attributed to the drag modeling, and not just numerical fluctuations between code executions.

The main takeaway of these results is that for LEO maneuvering, in timeframes of a few revolutions, the modeling of J2 effects is very important, but the modeling of drag (and supposedly, other even smaller perturbations) is not relevant. 

\begin{table}[htbp]\small
    \centering
    \begin{tabular}{clcccccc} \toprule
        \multicolumn{2}{c}{Scenario}& \multicolumn{2}{c}{Keplerian} & \multicolumn{2}{c}{J2}       & \multicolumn{2}{c}{J2+Drag} \\
                                 Name&$t_f$ (s)& \(n_i\) & \(\Delta v\) (m/s)  & \(n_i\) & \(\Delta v\) (m/s) & \(n_i\) & \(\Delta v\) (m/s)\\ \midrule
        Circle to circle rendez-vous         &3560.5& 2 & 887.56199 & 3 & 893.05336 & 3 & 893.05339 \\
        Noncoplanar rendez-vous  &11107.2& 4 & 36.14596  & 3 & 56.00653  & 3 & 56.01418 \\ \bottomrule
    \end{tabular}
    \caption{Comparison of results under different orbital models.}
    \label{tab:model_comparison}
\end{table}

\section{Conclusions}\label{sec:conc}

In summary, this work set out to optimize impulsive maneuvers in LEO, apply primer vector theory to the optimized maneuvers and adapt this theory to the types of relevant perturbations encountered in this space environment around Earth. The main theoretical contributions of this work were to derive primer vector optimality conditions through the modern theory of impulsive optimal control, extend primer vector theory to perturbed models and propose primer vector calculation methods that are valid under different perturbation classes: Glandorf, STM and ODE.

Numerically, it was found that the partition of coast times for the initial trajectory guess was a practical low-dimensional method of resampling initial guesses to find better local optima. It was also found that the CasADi model size heavily depends on the model, particularly the drag model used, as well as integrator choice. Balancing precision and weight is necessary with this method.

The orbital models can also be compared, with important differences being highlighted for orbital maneuver synthesis. Keplerian maneuvers are piecewise planar, differently to J2 trajectories, whose coasting arcs twist in space; these models have very different optima. Impulse instants, magnitudes and directions as well as their number can change under the J2 model. 

Despite the theoretical prediction that the STM method should not be applicable to the J2+Drag model, it was found that orbits in usual LEO heights have sufficiently low drag for the STM method to also be appliable. Slight differences in cost with respect to the J2 model do not justify the extra cost, and drag can be safely neglected in short time frames and usual LEO heights.

Future works can extend the time horizon for more complex maneuver cases or investigate aerobraking scenarios where certainly drag would have significant effects. The study of Very Low Earth Orbit Scenarios (VLEO) could complement the work developed here for cases with drag.

\textbf{Acknowledgements.}\quad The authors express their gratitude to E2MoC-Lab and ITA Space Center for all provided structures and facilities. This study was financed in part by the Coordenação de Aperfeiçoamento de Pessoal de Nível Superior - Brasil (CAPES) – Finance Code 001.

\begin{appendices}

\section{Uniform sampling of numbers with constant sum}\label{app:unif_sampl}

Sampling coast durations with constant total time is here simplified to a more abstract sampling problem with straightforward application. The goal is to generate a set of \(n\) numbers \((\x)_{1 \leq i \leq n} \in [0, 1]\) such that
\begin{equation}
    \sum_{i=1}^{n} \x_i = 1
\end{equation}
with uniform distribution on the set \(\Omega = \{\x \in [0, 1]^n \mid \sum_{i=1}^{n} \x_i = 1 \}\). The procedure relies on rejection sampling in \(n-1\) dimensions since the set \(\Omega\) has probability 0 in the space \([0, 1]^n\). Generate a uniform random sample \(\mathbf{z} \in [0, 1]^{n-1}\) and reject it if
\begin{equation}
    \sum_{i=1}^{n-1} \mathbf{z}_i > 1.
\end{equation}

Otherwise, accept the sample and set
\begin{align}
    \x_i &= \mathbf{z}_i, i = 1,\dots,n-1 \\
    \x_n &= 1 - \sum_{i=1}^{n-1} \x_i.
\end{align}

\section{Variable normalization}\label{app:var_norm}

Variable normalization was done by rescaling physical units for length and time before assembling the multiple shooting problem. Let \(L\) be a normalizing length, and \(T\) be a normalizing time. Then, a change of units is performed and new, rescaled physical parameters are found. Let \(\tilde{\bullet}\) represent the rescaled version of any variable \(\bullet\). This procedure is exemplified with the Keplerian model, but is analogous for the other models explored in the work:
\begin{equation}
    \ddot{\pos} = -\frac{\mu}{r^3} \pos \iff \frac{L}{T^2}\frac{d^2 \tilde{\pos}}{d \tilde{t}} = -\frac{\mu}{L^2 \tilde{r}^3} \tilde{\pos} \iff \frac{d^2 \tilde{\pos}}{d \tilde{t}} = -\frac{\tilde{\mu}}{\tilde{r}^3} \tilde{\pos}
\end{equation}
with \(\tilde{\mu} = \mu \frac{T^2}{L^3}\). The length scaling factor was taken to be the average of semimajor axes of the initial and final orbits, and the time scaling factor was set to the maneuver time.

\section{Smooth atmospheric model}\label{app:smooth_atm}

The US Standard Atmosphere model is defined piecewise~\cite{curtis_orbital_2020}, which is not compatible with Ipopt. The piecewise definition is based on the exponential interpolation heights $h_i$, where the density is $\rho_i$, with exponential coefficient $H_i$ for $i=1,\dots,28$. A smooth reformulation is therefore needed. This will be achieved in two steps.

Firstly, define the unit interval indicator function \(\mathbb{I} (x)\) as
\begin{equation}
    \mathbb{I} (x) = \begin{cases}
        1, 0 \leq x < 1 \\
        0, \text{ otherwise}
    \end{cases}
\end{equation}
and rewrite the piecewise definition given in~\cite{curtis_orbital_2020} for \(\rho(r)\) as
\begin{equation}
    \rho(r) = \sum_{i = 1}^{28} \rho_i \exp{\big(-\frac{\left(r - (h_i + R)\right)}{H_i}\big)} \mathbb{I}(\frac{r - R - h_i}{h_{i+1} - h_i}).
\end{equation}

Secondly, define a smooth indicator function \(\sigma_k(x)\) as
\begin{equation}
    \sigma_k(x) = \frac{1}{2} \left(\tanh{k x} + \tanh{k (1 - x)}\right)
\end{equation}
for some (large) value of \(k > 0\). Then, the smooth density model \(\rho_s(r)\) can be written
\begin{equation}
    \rho_s(r) = \sum_{i = 1}^{28} \rho_i \exp{\big(-\frac{\left(r - (h_i + R)\right)}{H_i}\big)} \sigma_k(\frac{r - R - h_i}{h_{i+1} - h_i}),
\end{equation}
thus giving an approximate, twice differentiable model as desired. With \(k = 200\), 
\begin{equation}
    \max_{r\in [R+100\text{ km}, R+1000\text{ km}]} \frac{\lvert \rho(r) - \rho_s(r) \rvert}{\rho(r)} \leq 2\%,
\end{equation}
which was deemed an acceptable approximation error.

\end{appendices}

\bibliography{sn-bibliography}

\end{document}